\documentclass[12pt]{article}
\usepackage{amsfonts}
\usepackage{mathrsfs}
\usepackage{amsmath,amssymb}

\def\a{\alpha}

\def\D{\Delta}

\def\LL{{\cal L}}

\def\Der{{\rm Der}}
\def\Im{{\rm Im}}
\def\Inn{{\rm Inn}}
\def\Ker{{\rm Ker}}
\def\v{\varphi}
\def\ssc{\scriptscriptstyle}
\def\cl{\centerline}
\def\rar{\rightarrow}

\def\vs{\vspace*}
\def\span{\mbox{Span}}
\def\even{{\bar0}}
\def\odd{{\bar1}}
\def\QQ{{\cal Q}}
\def\GG{{\cal G}}
\def\TT{\mathcal {T}}
\def\HH{\mathcal {H}}
\def\XX{\mathcal {X}}

\def\ZZ{\mathcal {Z}}
\def\ni{\noindent}
\def\Z{\mathbb{Z}{\ssc\,}}

\def\C{\mathbb{C}{\ssc\,}}
\def\b{\beta}
\def\g{\gamma}

\numberwithin{equation}{section}
\newtheorem{theo}{Theorem}[section]

\newtheorem{lemm}[theo]{Lemma}
\newtheorem{prop}[theo]{Proposition}
\newtheorem{clai}{Claim}

\def\adddot{$\!\!\!${\bf.}\ \ }

\def\SUM#1#2#3#4{\mbox{$\sum\limits#1{#2}#3{#4}$}}

\begin{document}

\cl{{\large \bf Topological $N\!=\!2$ superconformal
superbialgebras}
\footnote {Supported by NSF grants 10671027 of China\\[2pt]
\indent\ \,$^{*}$Corresponding author: jhZhou@seu.edu.cn}}\vs{6pt}

\cl{Lifang Lin$^{\dag)}$, Huanxia Fa$^{\ddag)}$, Jianhua
Zhou$^{*,\dag)}$}

\cl{\small $^{\dag)}$Department of Mathematics, Southeast
University, Nanjing 210096, China}

\cl{\small $^{\ddag)}$Department of Mathematics, Changshu Institute
of Technology, Changshu 215500, China}

\vs{6pt}

{\small
\parskip .005 truein
\baselineskip 3pt \lineskip 3pt

\noindent{{\bf Abstract.} Lie superbialgebra structures on the
centerless topological $N\!=\!2$ superconformal algebra $\TT$ are
considered, all of which are proved to be coboundary triangular.
\vs{5pt}

}

}

\vs{12pt}

 \cl{\bf\S1. \
Introduction}\setcounter{section}{1}\setcounter{equation}{0} The
notion of Lie bialgebras was first introduced by Drinfeld in 1983 in
a connection with Hamiltonian mechanics and Possion Lie group. Lie
bialgebras arising from non-abelian two dimensional Lie subalgebras
of Lie algebras (such as Witt, one-sided Witt, and Virasoro
algebras) have been considered in \cite{NT,T}. It is well known that
the $N=2$ superconformal algebras were first independently
constructed by Kac (see \cite{K}) and by Ademollo et al., which play
important roles in both mathematics and physics. The determination
of Lie bialgebra structures on a Lie algebra is considered to be the
first partial step towards quantizing Lie algebras. Lie
superbialgebra structures on the Ramond and twisted $N\!=\!2$
superconformal algebra were determined in \cite{YS} and \cite{FL}
respectively. In the present paper, we study Lie super-bialgebra
structures on the centerless topological $N\!=\!2$ superconformal
algebra. It is proved that all such Lie superbialgebras are
coboundary triangular.

Firstly we shall recall some related concepts. Let
$\TT=\TT^\even\oplus \TT^\odd$ be a super-vector space over the
complex number field $\C$. All elements in the following are assumed
to be $\Z_2$-homogeneous. For $x\in\TT$, the notation $[x]$ means
$x\in\TT^{[x]}$. Denote by $\tau$ and $\xi$ the {\it super-twist
map} of $\LL \otimes \LL$ and {\it super-cyclic map} of $\LL \otimes
\LL \otimes \LL$, i.e., $\tau(x_1\otimes x_2)= (-1)^{[x_1][x_2]}x_2
\otimes x_1$, $\xi=(1\otimes\tau)\cdot(\tau\otimes1): \,x_{1}
\otimes x_{2} \otimes x_{3}\mapsto (-1)^{[x_1]([x_2]+[x_3])}x_{2}
\otimes x_{3} \otimes x_{1}$ for $x_1,x_2,x_3\in \LL$.

A {\it Lie superalgebra} is a pair $(\TT,\v)$ consisting of a
super-vector space $\TT$ and a bilinear map $\v :\TT \otimes\TT
\to\TT$ such that $\v(\LL^i,\LL^j)\subset \LL^{i+j}$,
$\Ker(1\otimes1-\tau) \subset \Ker\,\v$, $\v \cdot (1 \otimes \v )
\cdot ({\bf 1}+\xi +\xi^{2})=0$. A {\it Lie super-coalgebra} is a
pair $(\TT,\D)$ consisting of a super-vector space $\TT$ and a
linear map $\D: \TT \rar \TT \otimes \TT$ such that
$\D(\LL^i)\!\subset\!\!\mbox{$\!\sum\limits_{j+k=i}$}\LL^{j}\otimes
\LL^k$, $\Im\,\D\subset\Im(1\otimes1- \tau)$, $({\bf 1}+\xi
+\xi^{2})\cdot(1\otimes\D)\cdot\D=0$. A {\it Lie super-bialgebra} is
a triple $(\TT, \v, \D )$ with $(\TT, \v)$ being a Lie superalgebra,
$(\TT,\D)$ a Lie super-coalgebra and $\D\v(x\otimes y)=x\cdot\D y-
(-1)^{[x][y]}y\cdot\D x$ for $x, y \in \TT$. The symbol ``$\cdot$''
means the adjoint diagonal action and $[x,y]=\v(x\otimes y)$, i.e.,
\begin{eqnarray*}
x\cdot(\SUM_{i}{a_{i}\otimes b_{i}})=\SUM_{i}({[x, a_{i}] \otimes
b_{i}+(-1)^{[x][a_i]}a_{i}\otimes[x,b_{i}]})\mbox{\ for }x,a_{i},
b_{i}\in\TT.
\end{eqnarray*}
A {\it coboundary Lie super-bialgebra} is a quadruple $(\TT , \v,
\D,r),$ where $(\TT , \v, \D)$ is a Lie super-bialgebra and
$\D=\D_r$ is a {\it coboundary of $r$} for $r\in\Im(1\otimes1-
\tau)\subset\TT\otimes\TT$. And $\D_r$ is defined by $\D_r(x)=
(-1)^{[r][x]}x\cdot r,\ \forall\,\,x \in\TT$. If $r=\sum_i{a_{i}
\otimes b_{i}}$, then $[r]=[a_i]+[b_i]$. A coboundary bialgebra is
called {\it triangular} if it satisfies the {\it classical
Yang-Baxter Equation}:
\begin{equation}
\label{e-CYBE} c(r): = [r^{12} , r^{13}] +[r^{12} , r^{23}] +[r^{13}
, r^{23}]=0,
\end{equation}
where $r=\sum_{i}{a_{i}\otimes b_{i}}$, $r^{12}=\sum_{i}a_{i}\otimes
b_{i} \otimes 1$, $r^{13}=\sum_{i}a_{i}\otimes 1\otimes b_{i}$ and
$r^{23}=\sum_{i}1\otimes a_{i}\otimes b_{i}$.

Let $V=V^{\bar0}\oplus V^{\bar1}$ be a $\TT$-module for the
superalgebra $\TT=\TT^{\bar0}\oplus\TT^{\bar1}$. A $\Z_2$-homogenous
linear map $d:\TT\to V$ is called a {\it homogenous derivations of
parity $[d]\in\Z_2$}, if $d(\LL^i)\subset V^{i+[d]}$ and there
exists $[d]$ such that $d([x,y])=(-1)^{[d][x]}x\cdot
d(y)-(-1)^{[y]([d]+[x])}y\cdot d(x)$ for $x,y\in \TT$. The
derivation $d$ is called {\it even} if $[d]=\bar0$, {\it odd} if
$[d]=\bar1$. Denote by $\Der^p(\TT,V)$ the set of homogenous
derivations of parity $p$. Let
$\Der(\TT,V)=\Der^{\bar0}(\TT,V)\oplus\Der^{\bar1}(\TT,V)$ be the
set of derivations from $\TT$ to $V$ and
$\Inn(\TT,V)=\Inn^{\bar0}(\TT,V)\oplus\Inn^{\bar1}(\TT,V)$ the set
of inner derivations, where $\Inn^p(\TT,V)$ is the set of {\it inner
derivations of parity $p$} consisting of $a_{\rm inn},$ $a\in V^p,$
defined by $a_{\rm inn}:x\mapsto (-1)^{[a][x]}x\cdot a$ for $x\in
\TT$, $[a]=p$.

Denote by $H^1(\TT,V)$ the {\it first cohomology group} of $\TT$
with coefficients in $V$. It is known
$H^1(\TT,V)\cong\Der(\TT,V)/\Inn(\TT,V)$. Say $r\in\TT$ satisfy the
{\it modified Yang-Baxter equation} if
\begin{eqnarray}
\label{e-MYBE} x\cdot c(r)=0,\ \ \forall\,\,x\in \TT.
\end{eqnarray}

We shall investigate the topological $N\!=\!2$ superconformal
algebra, $\TT\!=\!\TT^{\bar 0}\oplus\TT^{\bar 1}$, generated by the
bosonic operators $\LL_n$, $\HH_n$ and the fermionic operators
$\GG_n$, $\QQ_n$ with the following non-vanishing brackets
($\TT^{\bar 0}=\span_{\C}\{\LL_n,\HH_n\,|\,n\in\Z\}$ and $
\TT^{\bar1}=\span_{\C}\{\GG_n, \QQ_n|\,n\in\Z\}$)
\begin{eqnarray*}
&&[\LL_m,\QQ_n]=-n\QQ_{m+n},\\
&&[\LL_m,\HH_n]=-n\HH_{m+n},\ \ \ \ \ \ \ \ \ \ [\HH_m,
\GG_n]=\GG_{m+n},\\
&&[\LL_m,\GG_n]=(m-n)\GG_{m+n},\ \ \ \ \ \,[\HH_m,\QQ_n]=-\QQ_{m+n},\\
&&[\LL_m,\LL_n]= (m-n)\LL_{m+n},\ \ \ \ \ [\GG_m,\QQ_n
]=2\LL_{m+n}-2n\HH_{m+n}.
\end{eqnarray*}
\begin{theo}\adddot\label{main}
Every Lie super-bialgebra structure on the the centerless
topological $N\!=\!2$ superconformal algebra $\TT$ is triangular
coboundary.
\end{theo}

\cl{\bf\S3. \ Proof of main results}\setcounter{section}{3}
\setcounter{theo}{0} \setcounter{equation}{0}

The following lemma can be obtained by employing similar techniques
as \cite{NT,YS}.
\begin{lemm}\adddot\label{lemm3.1}
(i)\ \ The triple $(\TT,[\cdot,\cdot],\D_r)$ is a super-bialgebra if and only if $r$ satisfies (\ref{e-MYBE}).\\
(ii)\ \ If $r\in\TT\otimes\TT$ such that $x\cdot r=0$,
$\forall\,\,x\in\TT$, then $r=0$.\\
(iii)\ \ An element $r \in \Im(1\otimes1 - \tau)$ satisfies
(\ref{e-CYBE}) if only if it satisfies (\ref{e-MYBE}).
\end{lemm}
\begin{prop}\adddot\label{prop3.2}
$\Der(\TT ,V)=\Inn(\TT ,V)$, where $V=\TT\otimes \TT.$
\end{prop}
\ni{\it Proof.~}\def\b{j}\def\g{k} Note that
$V\!=\!\oplus_{i\in\Z}V_i$ is $\Z$-graded with $V_i=\sum_{j+k=i}
\TT_j\otimes\TT_k$. We say a derivation $d\in\Der(\TT ,V)$ is {\it
homogeneous of degree $i\in\Z$} if $d(V_j) \subset V_{i +j}$,
$\forall\,\,j \in\Z$ and denote such $i$ by $\deg d$. Set $\Der(\TT
, V)_i =\{d\in \Der(\TT , V) \,|\,{\rm deg\,}d =i\}$. For
$d\in\Der(\TT,V)$ and any $u\in\TT _j$, write
$d(u)=\sum_{k\in\Z}v_k\in V$ and set $d_i(v_j)=v_{i+j}$. Then
$d_i\in\Der(\TT ,V)_i$ and
\begin{eqnarray}\label{summable}
\mbox{$d = \sum\limits_{i \in\Z}$}d_i,\ \ d_i \in \Der(\TT, V)_i.
\end{eqnarray}
Such a sum is {\it summable}, i.e., for any $u \in\TT $ and $d(u) =
\sum_{i \in\Z} d_i(u)$, there are only finitely many $d_i(u)\neq 0$.
All the sums appear in the following are supposed to be summable.

\begin{clai}\adddot\label{clai1}
$d_i\in\Inn(\TT,V)$ and $d_0(L_0)=0$, $\forall\,\,i\in\Z^*$.
\end{clai}

For any $i\in\Z^*$, denote $u=-\frac1i{\ssc\,} d_{i}(L_0)\in V_{i}$.
For any $x_{j}\in \TT_{j}$, applying $d_{i}$ to $[L_0,x_{j}]=-j
x_{j},$ using $d_{i}(x_j)\in V_{i+j}$ and the action of $L_0$ on
$V_{i+j}$ is the scalar $-i-j$, we have
\begin{eqnarray}\label{equa-add-1}
-(i+j)d_{i}(x_{j}) - (-1)^{[d_{j}][x_{j}]}x_{j}\cdot d_{i}(L_0)=-j
d_{i}(x_{j}),
\end{eqnarray}
 i.e.,
$d_{i}(x_{j})=u_{\rm inn}(x_{j})$. Thus $d_{i}=u_{\rm inn}$ is
inner. Using (\ref{equa-add-1}) with $i=0$, we obtain $x\cdot
d_0(L_0)=0$. Thus by Lemma \ref{lemm3.1}, one has $d_0(L_0)=0$.
\begin{clai}\adddot\label{sub2}
By replacing $d_0$ by $d_0-u_{\rm inn}$ for some $u\in V_0$, one can
suppose $d_0=0$.
\end{clai}
For any $n\in\Z^*$, one can write
\begin{eqnarray*}
d_0(\LL_{n})\!\!\!&=\!\!\!&
\mbox{$\sum\limits_{i\in\Z}$}(a_{n,i}\LL_{i+n}\otimes\HH_{-i}+
a^\dag_{n,i}\HH_{i+n}\otimes\LL_{-i}+b_{n,i}\LL_{i+n}\otimes\GG_{-i}+
b^\dag_{n,i}\GG_{i+n}\otimes\LL_{-i}\\
&&+c_{n,i}\LL_{i+n}\otimes\QQ_{-i}+c^\dag_{n,i}\QQ_{i+n}\otimes\LL_{-i}
+d_{n,i}\HH_{i+n}\otimes\GG_{-i}+d^\dag_{n,i}\GG_{i+n}\otimes\HH_{-i}\\
&&+e_{n,i}\HH_{i+n}\otimes\QQ_{-i}+e^\dag_{n,i}\QQ_{i+n}\otimes\HH_{-i}
+f_{n,i}\GG_{i+n}\otimes\QQ_{-i}+f^\dag_{n,i}\QQ_{i+n}\otimes\GG_{-i}\\
&&+\a_{n,i}\LL_{i+n}\otimes\LL_{-i}+\beta_{n,i}\HH_{i+n}\otimes\HH_{-i}+\gamma_{n,i}\GG_{i+n}\otimes\GG_{-i}
+\mu_{n,i}\QQ_{i+n}\otimes\QQ_{-i}).
\end{eqnarray*}
Using the following facts,
$\LL_1\cdot(\XX_{i}\otimes\ZZ_{-i})=i\XX_{i}\otimes\ZZ_{1-i}-i\XX_{i+1}\otimes
\ZZ_{-i}\ \mbox{for}\ \XX,\ZZ=\HH\ \mbox{or}\ \QQ$,
\begin{eqnarray*}
&&\LL_1\cdot(\XX_{i}\otimes\ZZ_{-i})=(1-i)\XX_{i+1}\otimes
\ZZ_{-i}+i\XX_{i}\otimes\ZZ_{1-i}\ \mbox{for}\ \XX,\ZZ=\LL\ \mbox{or}\ \GG,\\
&&\LL_1\cdot(\XX_{i}\otimes\ZZ_{-i})=(1-i)\XX_{i+1}\otimes
\ZZ_{-i}+i\XX_{i}\otimes\ZZ_{1-i}\ \mbox{for}\ \XX=\LL\ \mbox{or}\
\GG,
\ \ZZ=\HH\ \mbox{or}\ \QQ,\\
&&\LL_1\cdot(\XX_{i}\otimes\ZZ_{-i})=(1+i)\XX_{i}\otimes\ZZ_{1-i}-i\XX_{i+1}\otimes
\ZZ_{-i}\ \mbox{for}\ \XX=\HH\ \mbox{or}\ \QQ,\ \ZZ=\LL\ \mbox{or}\
\GG,
\end{eqnarray*}
by subtracting a proper linear combination of the above inner
derivations, one can rewrite
$d_0(\LL_{1})-\mu_{1,-1}\QQ_{0}\otimes\QQ_{1}-\mu_{1,0}\QQ_{1}\otimes\QQ_{0}$
as
\begin{eqnarray*}
&&\!\!\!\!\!\!\!\!a_{1,-1}\LL_{0}\otimes\HH_{1}\!+\!a_{1,1}\LL_{2}\otimes\HH_{-1}\!+\!
a^\dag_{1,-2}\HH_{-1}\otimes\LL_{2}\!+\!a^\dag_{1,0}\HH_{1}\otimes\LL_{0}
\!+\!b_{1,-2}\LL_{-1}\otimes\GG_{2}\!+\!b_{1,1}\LL_{2}\otimes\GG_{-1}\\
&&\!\!\!\!\!\!\!\!+b^\dag_{1,-2}\GG_{-1}\otimes\LL_{2}\!+\!b^\dag_{1,1}\GG_{2}\otimes\LL_{-1}
\!+\!c_{1,-1}\LL_{0}\otimes\QQ_{1}\!+\!c_{1,1}\LL_{2}\otimes\QQ_{-1}
\!+\!c^\dag_{1,-2}\QQ_{-1}\otimes\LL_{2}\!+\!c^\dag_{1,0}\QQ_{1}\otimes\LL_{0}\\
&&\!\!\!\!\!\!\!\!+d_{1,-2}\HH_{-1}\otimes\GG_{2}\!+\!d_{1,0}\HH_{1}\otimes\GG_{0}
\!+\!d^\dag_{1,-1}\GG_{0}\otimes\HH_{1}
\!+\!d^\dag_{1,1}\GG_{2}\otimes\HH_{-1}
\!+\!e_{1,-1}\HH_{0}\otimes\QQ_{1}\!+\!e_{1,0}\HH_{1}\otimes\QQ_{0}\\
&&\!\!\!\!\!\!\!\!+e^\dag_{1,-1}\QQ_{0}\otimes\HH_{1}\!+\!e^\dag_{1,0}\QQ_{1}\otimes\HH_{0}
\!+\!f_{1,-1}\GG_{0}\otimes\QQ_{1}\!+\!f_{1,1}\GG_{2}\otimes\QQ_{-1}
\!+\!f^\dag_{1,-2}\QQ_{-1}\otimes\GG_{2}\!+\!f^\dag_{1,0}\QQ_{1}\otimes\GG_{0}\\
&&\!\!\!\!\!\!\!\!+\a_{1,-2}\LL_{-1}\otimes\LL_{2}\!+\!\a_{1,1}\LL_{2}\otimes\LL_{-1}
\!+\!\b_{1,-1}\HH_{0}\otimes\HH_{1}\!+\!\b_{1,0}\HH_{1}\otimes\HH_{0}
\!+\!\g_{1,-2}\GG_{-1}\otimes\GG_{2}\!+\!\g_{1,1}\GG_{2}\otimes\GG_{-1}.
\end{eqnarray*}
Applying $d_0$ to $[\LL_1,\LL_{-1}]= 2\LL_0$, we have
\begin{eqnarray}\label{81027n01}
\LL_1\cdot d_0(\LL_{-1})\!\!\!&=\!\!\!&\LL_{-1}\cdot d_0(\LL_1).
\end{eqnarray}
Comparing the coefficients of $\LL_{i}\otimes\HH_{-i}$ in
(\ref{81027n01}), we obtain $(i-2)a_{-1,i}=(i+1)a_{-1,i+1}$ for
$i\neq\pm1,0,2$, which together with the fact that the set $\{i\mid
a_{-1,i}\neq0\}$ is of finite rank, implies
\begin{eqnarray*}
a_{-1,i}=0\ \ \mbox{for}\ \ i\neq0,1,2\ \ \mbox{and}\ \
a_{-1,1}=-2a_{-1,0}=-2a_{-1,2},\ a_{1,-1}=a_{1,1}=0.
\end{eqnarray*}
Applying the same techniques to $\LL_{i}\otimes\QQ_{-i}$,
$\GG_{i}\otimes\HH_{-i}$ and $\GG_{i}\otimes\QQ_{-i}$, one has
\begin{eqnarray*}
x_{-1,i}=0,\ i\neq0,1,2,\ x_{-1,1}=-2x_{-1,0}=-2x_{-1,2},\
x_{1,-1}=x_{1,1}=0\ \ \mbox{for}\ \ x=c,f,d^\dag.
\end{eqnarray*}
Comparing the coefficients of $\HH_{i}\otimes\LL_{-i}$ in
(\ref{81027n01}), one has $(i-1)a^\dag_{-1,i}=(i+2)a^\dag_{-1,i+1}$
for $i\neq-2,0,\pm1$, which implies $a^\dag_{-1,i}=0$ for
$i\neq0,\pm1$, $2a^\dag_{-1,-1}=2a^\dag_{-1,1}=-a^\dag_{-1,0}$ and
$a^\dag_{1,0}=a^\dag_{1,-2}=0$. Comparing the coefficients of
$\QQ_{i}\otimes\LL_{-i}$, $\HH_{i}\otimes\GG_{-i}$ and
$\QQ_{i}\otimes\GG_{-i}$ in (\ref{81027n01}), one has
\begin{eqnarray*}
y_{-1,i}=0,\ i\neq0,\pm1\ \ \mbox{and}\ \
2y_{-1,-1}=2y_{-1,1}=-y_{-1,0},\ \ y_{1,0}=y_{1,-2}=0\ \ \mbox{for}\
\ y=d,c^\dag,f^\dag.
\end{eqnarray*}
Comparing the coefficients of $\LL_{i}\otimes\GG_{-i}$ and
$\GG_{i}\otimes\LL_{-i}$, we obtain
\begin{eqnarray*}
z_{-1,i}=z_{-1,1}+z_{-1,0}=
3z_{-1,-1}+z_{-1,0}+3z_{1,-2}=3z_{-1,2}-z_{-1,0}+3z_{1,1}=0,
\end{eqnarray*}
for $i\neq\pm1,0,2$ and $z=b$ or $z=b^\dag$. Finally, comparing the
coefficients of $\HH_{i}\otimes\QQ_{-i}$, $\QQ_{i}\otimes\HH_{-i}$,
$\LL_{i}\otimes\LL_{-i}$, $\GG_{i}\otimes\GG_{-i}$,
$\HH_{i}\otimes\HH_{-i}$ and $\QQ_{i}\otimes\QQ_{-i}$, we obtain
$x_{-1,i}=x_{-1,1}+x_{-1,0}+x_{1,-1}+x_{1,0}=0$,
\begin{eqnarray*}
&&y_{-1,j}=y_{-1,1}+y_{-1,0}=3y_{-1,-1}+y_{-1,0}+3y_{1,-2}=
3y_{-1,2}-y_{-1,0}+3y_{1,1}=0,
\end{eqnarray*}
for $i\neq0,1$ and $x=e,\beta,e^\dag,\mu$, $j\neq\pm1,0,2$ and
$y=\a,\gamma$. Then $d_0(\LL_{1})$ can be rewritten as
\begin{eqnarray*}
&&\!\!\!\!b_{1,-2}\LL_{-1}\!\otimes\!\GG_{2}\!+\!b_{1,1}\LL_{2}\!\otimes\!\GG_{-1}
\!+\!b^\dag_{1,-2}\GG_{-1}\!\otimes\!\LL_{2}\!+\!b^\dag_{1,1}\GG_{2}\!\otimes\!\LL_{-1}
\!+\!e_{1,-1}\HH_{0}\!\otimes\!\QQ_{1}\!+\!e_{1,0}\HH_{1}\!\otimes\!\QQ_{0}\\
&&\!\!\!\!+e^\dag_{1,-1}\QQ_{0}\otimes\HH_{1}+e^\dag_{1,0}\QQ_{1}\otimes\HH_{0}
+\a_{1,-2}\LL_{-1}\otimes\LL_{2}+\a_{1,1}\LL_{2}\otimes\LL_{-1}
+\beta_{1,-1}\HH_{0}\otimes\HH_{1}\\
&&\!\!\!\!+\beta_{1,0}\HH_{1}\otimes\HH_{0}+\gamma_{1,-2}\GG_{-1}\otimes\GG_{2}+\gamma_{1,1}\GG_{2}\otimes\GG_{-1}
+\mu_{1,-1}\QQ_{0}\otimes\QQ_{1}+\mu_{1,0}\QQ_{1}\otimes\QQ_{0}.
\end{eqnarray*}
For convenience, we introduce the following notations,
$u_{1}=\LL_{1}\otimes \HH_{-1}-\LL_0\otimes \HH_{0}$,
$u_{2}=\HH_{-1}\otimes \LL_1-\HH_0\otimes \LL_{0}$,
$u_{3}=\LL_{1}\otimes \QQ_{-1}-\LL_0\otimes \QQ_{0}$,
$u_{4}=\QQ_{-1}\otimes\LL_1-\QQ_0\otimes \LL_{0}$,
$u_{5}=\GG_{1}\otimes \HH_{-1}-\GG_0\otimes \HH_{0}$,
$u_{6}=\HH_{-1}\otimes \GG_1-\HH_0\otimes \GG_{0}$,
$u_{7}=\GG_{1}\otimes \QQ_{-1}-\GG_0\otimes \QQ_{0}$,
$u_{8}=\QQ_{-1}\otimes \GG_1-\QQ_0\otimes \GG_{0}$. Observing
\begin{eqnarray*}
&&\!\!\!\!\!\!\LL_{-1}\cdot u_{1}
=\LL_{-1}\!\otimes\!\HH_{0}-2\LL_{0}\!\otimes\!\HH_{-1}\!+\!\LL_{1}\!\otimes\!\HH_{-2},\
\LL_{-1}\cdot
u_{2}=-2\HH_{-1}\!\otimes\!\LL_{0}\!+\!\HH_{-2}\!\otimes\!\LL_{1}\!+\!\HH_{0}\!\otimes\!\LL_{-1},\\
&&\!\!\!\!\!\!\LL_{-1}\cdot
u_{3}=\LL_{-1}\!\otimes\!\QQ_{0}-2\LL_{0}\!\otimes\!\QQ_{-1}\!+\!\LL_{1}\!\otimes\!\QQ_{-2},\
\LL_{-1}\cdot
u_{4}=-2\QQ_{-1}\!\otimes\!\LL_{0}\!+\!\QQ_{-2}\!\otimes\!\LL_{1}\!+\!\QQ_{0}\!\otimes\!\LL_{-1},\\
&&\!\!\!\!\!\!\LL_{-1}\cdot
u_{5}=\GG_{-1}\!\otimes\!\HH_{0}-2\GG_{0}\!\otimes\!\HH_{-1}\!+\!\GG_{1}\!\otimes\!\HH_{-2},\
\LL_{-1}\cdot
u_{6}=-2\HH_{-1}\!\otimes\!\GG_{0}\!+\!\HH_{-2}\!\otimes\!\GG_{1}\!+\!\HH_{0}\!\otimes\!\GG_{-1},\\
&&\!\!\!\!\!\!\LL_{-1}\cdot
u_{7}=\GG_{-1}\!\otimes\!\QQ_{0}-2\GG_{0}\!\otimes\!\QQ_{-1}\!+\!\GG_{1}\!\otimes\!\QQ_{-2},\
\LL_{-1}\cdot
u_{8}=-2\QQ_{-1}\!\otimes\!\GG_{0}+\QQ_{-2}\!\otimes\!\GG_{1}+\QQ_{0}\!\otimes\!\GG_{-1},
\end{eqnarray*}
and $\LL_{1}\cdot u_{i}=\LL_{-1}\cdot u_{i}=0$, $\forall\,\,
i=1,\cdots,8$, one can replace $d_0$ by $d_0-\vartheta_{\rm inn}$
(where $\vartheta$ is a proper linear combination of $u_{i}$,
$i=1,\cdots,8$) and rewrite $d_0(\LL_{-1})$ as
\begin{eqnarray*}
&&\!\!\!\!\!\!\!\!b_{-1,0}\LL_{-1}\!\otimes\!\GG_{0}\!-\!b_{-1,0}\LL_{0}\!\otimes\!\GG_{-1}
\!-\!(b_{-1,0}/3+b_{1,-2})\LL_{-2}\!\otimes\!\GG_{1}\!+\!(b_{-1,0}/3\!-\!b_{1,1})\LL_{1}\!\otimes\!\GG_{-2}\\
&&\!\!\!\!\!\!\!\!+\gamma_{-1,0}\GG_{-1}\!\otimes\!\GG_{0}\!-\!(\gamma_{-1,0}/3+\gamma_{1,-2})\GG_{-2}\!\otimes\!\GG_{1}
\!-\!\gamma_{-1,0}\GG_{0}\!\otimes\!\GG_{-1}+(\gamma_{-1,0}/3\!-\!\gamma_{1,1})\GG_{1}\!\otimes\!\GG_{-2}\\
&&\!\!\!\!\!\!\!\!+b^\dag_{-1,0}\GG_{-1}\!\otimes\!\LL_{0}\!-\!b^\dag_{-1,0}\GG_{0}\!\otimes\!\LL_{-1}
\!-\!(b^\dag_{-1,0}/3+b^\dag_{1,-2})\GG_{-2}\!\otimes\!\LL_{1}
+(b^\dag_{-1,0}/3\!-\!b^\dag_{1,1})\GG_{1}\!\otimes\!\LL_{-2}\\
&&\!\!\!\!\!\!\!\!+\a_{-1,0}\LL_{-1}\!\otimes\!\LL_{0}\!-\!(\a_{-1,0}/3+\a_{1,-2})\LL_{-2}\!\otimes\!\LL_{1}
\!-\!\a_{-1,0}\LL_{0}\!\otimes\!\LL_{-1}+(\a_{-1,0}/3\!-\!\a_{1,1})\LL_{1}\!\otimes\!\LL_{-2}\\
&&\!\!\!\!\!\!\!\!+e_{-1,0}\HH_{-1}\!\!\otimes\!\!\QQ_{0}\!\!-\!\!(e_{-1,0}+e_{1,-1}+e_{1,0})\HH_{0}\!\otimes\!\QQ_{-1}
\!\!+\!\!e^\dag_{-1,0}\QQ_{-1}\!\otimes\!\HH_{0}\!\!-\!\!(e^\dag_{-1,0}\!+\!e^\dag_{1,-1}+e^\dag_{1,0})\QQ_{0}\!\otimes\!\HH_{-1}\\
&&\!\!\!\!\!\!\!\!+\beta_{-1,0}\HH_{-1}\!\otimes\!\HH_{0}\!-\!(\beta_{-1,0}\!\!+\!\!\beta_{1,-1}+\beta_{1,0})\HH_{0}\!\otimes\!\HH_{-1}
\!+\!\mu_{-1,0}\QQ_{-1}\!\otimes\!\QQ_{0}\!\!-\!\!(\mu_{-1,0}+\mu_{1,-1}\!\!+\!\!\mu_{1,0})\QQ_{0}\!\otimes\!\QQ_{-1}.
\end{eqnarray*}
Applying $d_0$ to $[\LL_1,\LL_{-2}]= 3\LL_{-1}$, one has
\begin{eqnarray}\label{81028n01}
\LL_1\cdot d_0(\LL_{-2})\!\!\!&=\!\!\!& \LL_{-2}\cdot
d_0(\LL_1)+3d_0(\LL_{-1}).
\end{eqnarray}
Comparing the coefficients of $\LL_{i-1}\otimes\HH_{-i}$ in
(\ref{81028n01}), one has $(i-3)a_{-2,i}=(i+1)a_{-2,i+1}$ for
$i\neq0,1,2$, which together with $\{i\mid\a_{-2,i}\neq0\}$ being
finite, implies
\begin{eqnarray*}
a_{-2,i}=a_{-2,1}+3a_{-2,0}=a_{-2,2}-3a_{-2,0}=a_{-2,3}+a_{-2,0}=0,\
\ \forall\,\,i\neq0,1,2,3.
\end{eqnarray*}
Similarly, comparing the coefficients of $\LL_{i-1}\otimes\QQ_{-i}$,
 $\GG_{i-1}\otimes\HH_{-i}$, $\GG_{i-1}\otimes\QQ_{-i}$, we obtain
\begin{eqnarray*}
x_{-2,i}=x_{-2,1}+3x_{-2,0}=x_{-2,2}-3x_{-2,0}=x_{-2,3}+x_{-2,0}=0,
\,\forall\,\,i\neq0,1,2,3,\,x=c,f,d^\dag.
\end{eqnarray*}
Comparing the coefficients of $\HH_{i-1}\otimes\LL_{-i}$ in
(\ref{81028n01}), one has $(i-2)a^\dag_{-2,i}=(i+2)a^\dag_{-2,i+1}$
for $i\neq\pm1,0$, which implies
$a^\dag_{-2,i}=a^\dag_{-2,1}+a^\dag_{-2,0}
=3a^\dag_{-2,-1}+a^\dag_{-2,0}=3a^\dag_{-2,2}-a^\dag_{-2,0}=0$,
$\forall\,\,i\neq\pm1,0,2$. Then comparing the coefficients of
$\QQ_{i-1}\otimes\LL_{-i}$, $\HH_{i-1}\otimes\GG_{-i}$,
$\QQ_{i-1}\otimes\GG_{-i}$, we obtain $y_{-2,i}=y_{-2,1}+y_{-2,0}
=3y_{-2,-1}+y_{-2,0}=3y_{-2,2}-y_{-2,0}=0$,
$\forall\,\,i\neq\pm1,0,2,\,y=d,c^\dag,f^\dag$. Comparing the
coefficients of $\LL_{i-1}\otimes\GG_{-i}$ and
$\GG_{i-1}\otimes\LL_{-i}$ in (\ref{81028n01}), one has
\begin{eqnarray*}
(i+2)z_{-2,i+1}=(i-3)z_{-2,i}\ \
\mbox{for}\,\,i\neq\pm2,\pm1,0,3,\,\,z=b,b^\dag,
\end{eqnarray*}
which implies $z_{-2,i}=2z_{-2,1}+3z_{-2,0}-3z_{-1,0}=
z_{-2,2}-z_{-2,0}+2z_{-1,0}=4z_{-2,-1}+z_{-2,0}+z_{-1,0}=
4z_{-2,3}+z_{-2,0}-3z_{-1,0}=z_{1,-2}=z_{1,1}=0$, for all
$i\neq\pm1,0,2,3$ and $z=b,b^\dag$. Comparing the coefficients of
$\HH_{i-1}\otimes\QQ_{-i}$ and $\QQ_{i-1}\otimes\HH_{-i}$ in
(\ref{81028n01}), one has $(i-2)e_{-2,i}=(i+1)e_{-2,i+1}$ for
$i\neq0,1,\ x=e,e^\dag$, which implies
$x_{-2,i}=x_{-2,1}+2x_{-2,0}-3x_{-1,0}+x_{1,0}=
x_{-2,2}-x_{-2,0}+3x_{-1,0}+2x_{1,-1}+x_{1,0}=0$ for all
$i\neq\pm1,0,2$ and $x=e,e^\dag$. Comparing the coefficients of
$\LL_{i-1}\otimes\LL_{-i}$, $\GG_{i-1}\otimes\GG_{-i}$,
$\HH_{i-1}\otimes\HH_{-i}$ and $\QQ_{i-1}\otimes\QQ_{-i}$ in
(\ref{81028n01}), we obtain
\begin{eqnarray*}
y_{-2,i}\!\!\!&=&\!\!\!2y_{-2,1}+3y_{-2,0}-3y_{-1,0}=
y_{-2,2}-y_{-2,0}+2y_{-1,0}\\
\!\!\!&=&\!\!\!4y_{-2,3}+y_{-2,0}-3y_{-1,0}=4y_{-2,-1}+y_{-2,0}+y_{-1,0}=y_{1,1}=y_{1,-2}=0,\\
z_{-2,j}\!\!\!&=&\!\!\!z_{-2,1}+2z_{-2,0}+z_{1,0}-3z_{-1,0}=
z_{-2,2}-z_{-2,0}+2z_{1,-1}+3z_{-1,0}+z_{1,0}=0,
\end{eqnarray*}
where $i\neq\pm1,0,2,3$, $j\neq0,1,2$, $y=\a,\gamma$ and
$z=\beta,\mu$. From $d_0([\LL_2,\LL_{-1}])=3d_0\LL_{1}$, one has
\begin{eqnarray}\label{1029n11}
\LL_{-1}\cdot d_0(\LL_2)\!\!\!&=\!\!\!&\LL_{2}\cdot
d_0(\LL_{-1})-3d_0(\LL_1).
\end{eqnarray}
Comparing the coefficients of $\LL_{i+1}\otimes\HH_{-i}$ in
(\ref{1029n11}), one has $(i+3)a_{2,i}=(i-1)a_{2,i-1}$ for
$i\neq\pm1,0,2$, which together with $\{i\mid\a_{2,i}\neq0\}$ being
finite, implies
\begin{eqnarray*}
a_{2,i}=a_{2,-1}+3a_{2,0}=a_{2,-2}-3a_{2,0}=a_{2,-3}+a_{2,0}=0,\
\forall\,\,i\neq-3,\cdots,0.
\end{eqnarray*}
Similarly, comparing the coefficients of $\LL_{i+1}\otimes\QQ_{-i}$,
$\GG_{i+1}\otimes\HH_{-i}$, $\GG_{i+1}\otimes\QQ_{-i}$, we obtain
\begin{eqnarray*}
x_{2,i}=x_{2,-1}+3x_{2,0}=x_{2,-2}-3x_{2,0}=x_{2,-3}+x_{2,0}=0,\
\forall\,\,i\neq-3,\cdots,0,\ x=c,f,d^\dag.
\end{eqnarray*}
Comparing the coefficients of $\HH_{i+1}\otimes\LL_{-i}$,
$\QQ_{i+1}\otimes\LL_{-i}$, $\HH_{i+1}\otimes\GG_{-i}$ and
$\QQ_{i+1}\otimes\GG_{-i}$, one has
\begin{eqnarray*}
y_{2,i}=3y_{2,1}+y_{2,0}=y_{2,-1}+y_{2,0}= 3y_{2,-2}-y_{2,0}=0,\
\forall\,\,i\neq0,\pm1,-2,\ y=a^\dag,c^\dag,d,f^\dag.
\end{eqnarray*}
Comparing the coefficients of $\LL_{i+1}\otimes\GG_{-i}$ and
$\GG_{i+1}\otimes\LL_{-i}$ in (\ref{1029n11}), we obtain
\begin{eqnarray*}
z_{2,i}=z_{-1,0}=2z_{2,-1}+3z_{2,0}=4z_{2,1}+z_{2,0}=4z_{2,-3}+z_{2,0}=0,\
z_{2,-2}=z_{2,0},
\end{eqnarray*}
for all $i\neq0,\pm1,-3$ and $z=b,b^\dag$. Comparing the
coefficients of $\HH_{i+1}\otimes\QQ_{-i}$ and
$\QQ_{i+1}\otimes\HH_{-i}$ in (\ref{1029n11}), one has
$e_{2,i}=e_{2,-1}+2e_{2,0}+e_{-1,0}-3e_{1,0}=
e_{2,-2}-e_{2,0}-e_{-1,0}-2e_{1,-1}+e_{1,0}=0$, for all
$i\neq0,-1,-2$ and $x=e,e^\dag$. Comparing the coefficients of
$\LL_{i+1}\otimes\LL_{-i}$ and $\GG_{i+1}\otimes\GG_{-i}$ in
(\ref{1029n11}), one has
$y_{2,i}=y_{-1,0}=2y_{2,-1}+3y_{2,0}=y_{2,-2}-y_{2,0}
=4y_{2,-3}+y_{2,0}=4y_{2,1}+y_{2,0}=0$, for all $i\neq-3,-2,0,\pm1$
and $y=\a,\gamma$. Comparing the coefficients of
$\HH_{i+1}\otimes\HH_{-i}$ and $\QQ_{i+1}\otimes\QQ_{-i}$ in
(\ref{1029n11}), one has
$z_{2,i}=z_{2,-1}+2z_{2,0}+z_{-1,0}-3z_{1,0}=
z_{2,-2}-z_{2,0}-z_{-1,0}-2z_{1,-1}+z_{1,0}=0$, for all
$i\neq-2,-1,0$ and $z=\beta,\mu$. Denote
$v_1=2\LL_{0}\otimes\GG_{0}-\LL_1\otimes\GG_{-1}-\LL_{-1}\otimes\GG_1$,
$v_2=2\GG_{0}\otimes\LL_{0}-\GG_1\otimes\LL_{-1}-\GG_{-1}\otimes\LL_1$,
$v_3=2\LL_{0}\otimes\LL_{0}-\LL_1\otimes\LL_{-1}-\LL_{-1}\otimes\LL_1$,
$v_4=2\GG_{0}\otimes\GG_{0}-\GG_1\otimes\GG_{-1}-\GG_{-1}\otimes\GG_1$.
Observing the facts, $\LL_{-2}\cdot
v_1=-4\LL_{-2}\otimes\GG_{0}+6\LL_{-1}\otimes\GG_{-1}-4\LL_{0}\otimes\GG_{-2}
+\LL_{-3}\otimes\GG_{1}+\LL_{1}\otimes\GG_{-3}$,
\begin{eqnarray*}
&&\LL_{-2}\cdot
v_2=-4\GG_{-2}\otimes\LL_{0}+6\GG_{-1}\otimes\LL_{-1}-4\GG_{0}\otimes\LL_{-2}
+\GG_{-3}\otimes\LL_{1}+\GG_{1}\otimes\LL_{-3},\\
&&\LL_{-2}\cdot
v_3=-4\LL_{-2}\otimes\LL_{0}+6\LL_{-1}\otimes\LL_{-1}-4\LL_{0}\otimes\LL_{-2}
+\LL_{-3}\otimes\LL_{1}+\LL_{1}\otimes\LL_{-3},\\
&&\LL_{-2}\cdot
v_4=-4\GG_{-2}\otimes\GG_{0}+6\GG_{-1}\otimes\GG_{-1}-4\GG_{0}\otimes\GG_{-2}
+\GG_{-3}\otimes\GG_{1}+\GG_{1}\otimes\GG_{-3}
\end{eqnarray*}
and $\LL_{1}\cdot v_i=\LL_{-1}\cdot v_i=0$,
$\forall\,\,i=1,\cdots,4$, one can replace $d_0$ by $d_0-\varpi_{\rm
inn}$ (where $\varpi$ is a proper linear combination of $v_i$
$(i=1,\cdots,4)$) and rewrite $d_0(\LL_{-2})$ as
\begin{eqnarray*}
&&\!\!\!\!\!\!e_{-2,0}\HH_{-2}\!\otimes\!\QQ_{0}\!-\!(2e_{-2,0}-3e_{-1,0}+e_{1,0})\HH_{-1}\!\otimes\!\QQ_{-1}
\!+\!(e_{-2,0}-3e_{-1,0}-2e_{1,-1}-e_{1,0})\HH_{0}\!\otimes\!\QQ_{-2}\\
&&\!\!\!\!\!\!+e^\dag_{-2,0}\QQ_{-2}\!\otimes\!\HH_{0}\!
-\!(2e^\dag_{-2,0}-3e^\dag_{-1,0}+e^\dag_{1,0})\QQ_{-1}\!\otimes\!\HH_{-1}
\!+(e^\dag_{-2,0}-3e^\dag_{-1,0}\!-\!2e^\dag_{1,-1}-e^\dag_{1,0})\QQ_{0}\!\otimes\!\HH_{-2}\\
&&\!\!\!\!\!\!+\beta_{-2,0}\HH_{-2}\!\otimes\!\HH_{0}\!-\!(2\beta_{-2,0}\!+\!\beta_{1,0}-3\beta_{-1,0})\HH_{-1}\!\otimes\!\HH_{-1}
\!+\!(\beta_{-2,0}\!-\!2\beta_{1,-1}\!-\!3\beta_{-1,0}\!-\!\beta_{1,0})\HH_{0}\!\otimes\!\HH_{-2}\!\\
&&\!\!\!\!\!\!+\mu_{-2,0}\QQ_{-2}\!\otimes\!\QQ_{0}\!-\!(2\mu_{-2,0}
\!\!+\!\!\mu_{1,0}\!\!-\!\!3\mu_{-1,0})\QQ_{-1}\!\otimes\!\QQ_{-1}\!+\!(\mu_{-2,0}\!-\!2\mu_{1,-1}-3\mu_{-1,0}-\mu_{1,0})\QQ_{0}\!\otimes\!\QQ_{-2}.
\end{eqnarray*}

Applying $d_0$ to $[\LL_2,\LL_{-2}]= 4\LL_0$, we have
\begin{eqnarray}\label{81030n01}
\LL_2\cdot d_0(\LL_{-2})\!\!\!&=\!\!\!& \LL_{-2}\cdot d_0(\LL_2).
\end{eqnarray}
Comparing the coefficients of $\LL_2\otimes\HH_{-2}$,
$\LL_{-2}\otimes\HH_2$, $\LL_2\otimes\QQ_{-2}$,
$\LL_{-2}\otimes\QQ_2$, $\GG_2\otimes\HH_{-2}$,
$\GG_{-2}\otimes\HH_2$, $\GG_2\otimes\QQ_{-2}$,
$\GG_{-2}\otimes\QQ_2$ in (\ref{81030n01}), one has
$a_{2,0}=a_{-2,0}=c_{2,0}=c_{-2,0}=d^\dag_{2,0}=d^\dag_{-2,0}=f_{2,0}=f_{-2,0}=0$.
Comparing the coefficients of $\HH_{-2}\otimes\LL_2$,
$\HH_2\otimes\LL_{-2}$, $\QQ_{-2}\otimes\LL_2$,
$\QQ_2\otimes\LL_{-2}$, $\HH_{-2}\otimes\GG_2$,
$\HH_2\otimes\GG_{-2}$, $\QQ_{-2}\otimes\GG_2$,
$\QQ_2\otimes\GG_{-2}$ in (\ref{81030n01}), one has
$a^\dag_{-2,0}=a^\dag_{2,0}=c^\dag_{-2,0}=c^\dag_{2,0}=d_{-2,0}
=d_{2,0}=f^\dag_{-2,0}=f^\dag_{2,0}=0$. Comparing the coefficients
of $\LL_0\otimes\GG_0$, $\GG_0\otimes\LL_0$, $\LL_0\otimes\LL_0$ and
$\GG_0\otimes\GG_0$ in (\ref{81030n01}), one has
$b_{-2,0}=b^\dag_{-2,0}=\a_{2,0}=\g_{2,0}=0$. Comparing the
coefficients of $\HH_1\otimes\QQ_{-1}$, $\QQ_1\otimes\HH_{-1}$,
$\HH_0\otimes\HH_0$ and $\QQ_0\otimes\QQ_0$ in (\ref{81030n01}), one
has
\begin{eqnarray}\label{81031n01}
x_{-2,0}+x_{2,0}=x_{-1,0}+x_{1,0},\
y_{-2,0}+y_{2,0}=y_{-1,0}+y_{1,0}\ \ \mbox{for}\,\,x=e,\beta,\
y=e^\dag,\mu.
\end{eqnarray}
Thus one can rewrite
$d_0(\LL_{1})-\!\mu_{1,-1}\QQ_{0}\!\otimes\!\QQ_{1}-\mu_{1,0}\QQ_{1}\!\otimes\!\QQ_{0}$
and $d_0(\LL_{-1})$ respectively as
\begin{eqnarray*}
&&\!\!\!\!\!\!\!\!e_{1,-1}\HH_{0}\!\otimes\!\QQ_{1}\!+\!e_{1,0}\HH_{1}\!\otimes\!\QQ_{0}
\!\!+e^\dag_{1,-1}\QQ_{0}\!\otimes\!\HH_{1}\!+\!e^\dag_{1,0}\QQ_{1}\!\otimes\!\HH_{0}
\!+\b_{1,-1}\HH_{0}\!\otimes\!\HH_{1}\!+\!\b_{1,0}\HH_{1}\!\otimes\!\HH_{0},\\
&&\!\!\!\!\!\!\!\!e_{-1,0}\HH_{-1}\!\otimes\!\QQ_{0}\!-\!(e_{-1,0}+e_{1,-1}+e_{1,0})\HH_{0}\!\otimes\!\QQ_{-1}
\!+\!e^\dag_{-1,0}\QQ_{-1}\!\otimes\!\HH_{0}\!-\!(e^\dag_{-1,0}+e^\dag_{1,-1}+e^\dag_{1,0})\QQ_{0}\!\otimes\!\HH_{-1}\\
&&\!\!\!\!\!\!\!\!+\!\beta_{-1,0}\HH_{-1}\!\otimes\!\HH_{0}\!-\!(\beta_{-1,0}\!+\!\beta_{1,-1}\!+\!\beta_{1,0})\HH_{0}\!\otimes\!\HH_{-1}
\!+\!\mu_{-1,0}\QQ_{-1}\!\otimes\!\!\QQ_{0}\!-\!\!(\mu_{-1,0}\!+\!\mu_{1,-1}\!+\!\mu_{1,0})\QQ_{0}\!\otimes\!\QQ_{\!-1\!}.
\end{eqnarray*}

One can write
\begin{eqnarray*}
d_0(\GG_{n})\!\!\!&=\!\!\!&
\mbox{$\sum\limits_{i\in\Z}$}(m_{n,i}\LL_{i+n}\otimes\HH_{-i}+
m^\dag_{n,i}\HH_{i+n}\otimes\LL_{-i}+k_{n,i}\LL_{i+n}\otimes\GG_{-i}+
k^\dag_{n,i}\GG_{i+n}\otimes\LL_{-i}\\
&&+u_{n,i}\LL_{i+n}\otimes\QQ_{-i}+u^\dag_{n,i}\QQ_{i+n}\otimes\LL_{-i}
+p_{n,i}\HH_{i+n}\otimes\GG_{-i}+p^\dag_{n,i}\GG_{i+n}\otimes\HH_{-i}\\
&&+q_{n,i}\HH_{i+n}\otimes\QQ_{-i}+q^\dag_{n,i}\QQ_{i+n}\otimes\HH_{-i}
+s_{n,i}\GG_{i+n}\otimes\QQ_{-i}+s^\dag_{n,i}\QQ_{i+n}\otimes\GG_{-i}\\
&&+t_{n,i}\LL_{i+n}\otimes\LL_{-i}+w_{n,i}\HH_{i+n}\otimes\HH_{-i}+t^\dag_{n,i}\GG_{i+n}\otimes\GG_{-i}
+w^\dag_{n,i}\QQ_{i+n}\otimes\QQ_{-i}).
\end{eqnarray*}
Applying $d_0$ to $[\LL_1,\GG_1]=0$, we have $\LL_1\cdot
d_0(\GG_{1})=(-1)^{[d]}\GG_1\cdot d_0(\LL_1)$. Comparing the
coefficients of $\LL_{i+1}\otimes\HH_{1-i}$, one has
$(i-1)m_{1,i-1}=im_{1,i}\ \ \mbox{for}\ i\neq0,1$, which together
with $\{i\,|\,m_{1,i}\neq0\}$ being finite, implies
$m_{1,i}=e^\dag_{1,-1}=e^\dag_{1,0}=0$, $\forall\,\,i\in\Z^*$. Using
the same techniques to $\LL_{i+1}\otimes\QQ_{1-i}$,
$\GG_{i+1}\otimes\HH_{1-i}$, $\GG_{i+1}\otimes\QQ_{1-i}$, one has
\begin{eqnarray*}
x_{1,i}=y_{1,-1}=y_{1,0}=0,\ \mbox{for}\,\,x=u,p^\dag,s,\
y=\mu,\beta,e,\  \forall\,\,i\in\Z^*.
\end{eqnarray*}
Comparing the coefficients of $\HH_{i+1}\otimes\LL_{1-i}$,
$\QQ_{i+1}\otimes\LL_{1-i}$, $\HH_{i+1}\otimes\GG_{1-i}$,
$\QQ_{i+1}\otimes\GG_{1-i}$, one has
$m^\dag_{1,i}=u^\dag_{1,i}=p_{1,i}=s^\dag_{1,i}=0,\
\forall\,\,i\,\,\neq -1$. Comparing the coefficients of
$\LL_{i+1}\otimes\GG_{1-i}$, $\GG_{i+1}\otimes\LL_{1-i}$,
$\LL_{i+1}\otimes\LL_{1-i}$, $\GG_{i+1}\otimes\GG_{1-i}$, one has
$x_{1,i}=0,\ x_{1,-1}=-x_{1,0}$ for $x=k,t,k^\dag,t^\dag$,
$\forall\,\,i\neq0,-1$. Comparing the coefficients of
$\HH_{i+1}\otimes\QQ_{1-i}$ , $\QQ_{i+1}\otimes\HH_{1-i}$,
$\HH_{i+1}\otimes\HH_{1-i}$, $\QQ_{i+1}\otimes\QQ_{1-i}$, one has
$q_{1,i}=q^\dag_{1,i}=w_{1,i}=w^\dag_{1,i}=0$ for $i\in\Z$. Also by
the above results, we have $d_0(\LL_{1})=0$.

Applying $d_0$ to $[\LL_1,\GG_0]=\GG_1$, we have $\LL_1\cdot
d_0(\GG_0)=d_0(\GG_1)$. Comparing the coefficients of
$\LL_{i+1}\otimes\HH_{-i}$, one has $(i-1)m_{0,i}=(i+1)m_{0,i+1}\ \
\mbox{for}\ i\neq0$. Since $\{i\mid\ m_{0,i}\neq0\}$ is a finite
set, we obtain $m_{0,i}=0,\ m_{0,1}=-m_{0,0}+m_{1,0}$ for
$i\neq0,1$. Applying the same techniques to
$\LL_{i+1}\otimes\QQ_{-i}$, $\GG_{i+1}\otimes\HH_{-i}$,
$\GG_{i+1}\otimes\QQ_{-i}$, $\HH_{i}\otimes\LL_{1-i}$,
$\QQ_{i}\otimes\LL_{1-i}$, $\HH_{i}\otimes\GG_{1-i}$,
$\QQ_{i}\otimes\GG_{1-i}$, we have $x_{0,i}=0,\
x_{0,1}=-x_{0,0}+x_{1,0}$ and $y_{0,j}=0,\
y_{0,-1}=-y_{0,0}+y_{1,-1}$ for
$x=u,p^\dag,s,y=m^\dag,u^\dag,p,s^\dag$,
$\forall\,\,i\neq0,1,j\neq0,-1$. Comparing the coefficients of
$\LL_{i+1}\otimes\GG_{-i}$, $\GG_{i+1}\otimes\LL_{-i}$,
$\LL_{i+1}\otimes\LL_{-i}$, $\GG_{i+1}\otimes\GG_{-i}$, we have
\begin{eqnarray*}
&&x_{0,i}=x_{0,1}+\frac{1}{2}x_{0,0}-\frac{1}{2}x_{1,0}
=x_{0,-1}+\frac{1}{2}x_{0,0}+\frac{1}{2}x_{1,0}=0\
\mbox{for}\,\,x=k,k^\dag,t,t^\dag,\ \forall\,\,i\neq0,-1,1.
\end{eqnarray*}
Comparing the coefficients of $\HH_{i}\otimes\QQ_{1-i}$,
$\QQ_{i}\otimes\HH_{1-i}$, $\HH_{i}\otimes\HH_{1-i}$,
$\QQ_{i}\otimes\QQ_{1-i}$, one has
\begin{eqnarray*}
q_{0,i}=q^\dag_{0,i}=w_{0,i}=w^\dag_{0,i}=0,\ \forall\,\,i\neq0.
\end{eqnarray*}

Applying $d_0$ to $[\LL_{-1},\GG_1]=-2\GG_0$, we have $\LL_{-1}\cdot
d_0(\GG_1)+2d_0(\GG_0)=(-1)^{[d]}\GG_1\cdot d_0(\LL_{-1})$.
Comparing the coefficients of $\LL_0\otimes\HH_0,$
$\HH_0\otimes\LL_0$, $\LL_{0}\otimes\QQ_{0}$,
$\QQ_{0}\otimes\LL_{0}$, one has
\begin{eqnarray*}
&&m_{1,0}=m_{0,0}\!\!-\!\!(-1)^{[d]}e^\dag_{-1,0},\
m^\dag_{1,-1}=m^\dag_{0,0}\!\!+\!\!(-1)^{[d]}e_{-1,0},\\
&&u_{1,0}=-u_{0,0}\!\!+\!\!(-1)^{[d]}\mu_{-1,0},\
u^\dag_{1,-1}=u^\dag_{0,0}\!\!+\!\!(-1)^{[d]}\mu_{-1,0}.
\end{eqnarray*}
Comparing the coefficients of $\HH_{0}\otimes\GG_{0}$,
$\GG_{0}\otimes\HH_{0}$, $\GG_{0}\otimes\QQ_{0}$,
$\QQ_{0}\otimes\GG_{0}$, $\LL_{0}\otimes\GG_{0}$,
$\GG_{0}\otimes\LL_{0}$, $\LL_{0}\otimes\LL_{0}$,
$\GG_{0}\otimes\GG_{0}$, we have
$y_{0,0}=y^\dag_{0,0}=t_{0,0}=t^\dag_{0,0}=0$,
$p^\dag_{1,0}=p^\dag_{0,0}\!\!+\!\!(-1)^{[d]}\frac{1}{2}\beta_{\!\!-1,0},\
\ p_{1,-1}=p_{0,0}\!\!+\!\!(-1)^{[d]+1}\frac{1}{2}\beta_{\!\!-1,0}$
and $s_{1,0}=s_{0,0}\!\!+\!\!(-1)^{[d]}\frac{1}{2}e_{\!\!-1,0},\
s^\dag_{1,-1}=s^\dag_{0,0}\!\!+\!\!(-1)^{[d]}\frac{1}{2}e^\dag_{\!\!-1,0}$.
Comparing the coefficients of $\HH_0\otimes\QQ_0$,
$\QQ_0\otimes\HH_0$, $\HH_{0}\otimes\HH_{0}$,
$\QQ_{0}\otimes\QQ_{0}$, one has
\begin{eqnarray*}
q^\dag_{0,0}=q_{0,0}=(-1)^{[d]}\mu_{-1,0}, \ \
w_{0,0}=(-1)^{[d]+1}2e_{-1,0}+(-1)^{[d]}2e^\dag_{-1,0},\ \
w^\dag_{0,0}=0.
\end{eqnarray*}
Applying $d_0$ to $[\LL_{-1},\GG_0]\!=\!-\GG_{-1}$, one has
$\LL_{-1}\cdot d_0(\GG_0)+d_0(\GG_{-1})\!=\!(-1)^{[d]}\GG_0\cdot
d_0(\LL_{-1})$. Comparing the coefficients of
$\LL_{-1}\otimes\HH_0$, $\LL_0\otimes\HH_{-1}$,
$\LL_1\otimes\HH_{-2}$, $\HH_0\otimes\LL_{-1}$,
$\HH_{-2}\otimes\LL_1$, $\HH_{-1}\otimes\LL_{0}$, we have
$m_{-1,0}=m_{0,0}+(-1)^{[d]}2e^\dag_{-1,0},\
m_{-1,1}=(-1)^{[d]+1}4e^\dag_{-1,0},\
m_{-1,2}=(-1)^{[d]}e^\dag_{-1,0}$,
\begin{eqnarray*}
&&m^\dag_{-1,1}=m^\dag_{0,0}+(-1)^{[d]+1}2e_{-1,0},\
m^\dag_{-1,-1}=(-1)^{[d]+1}e_{-1,0},\
m^\dag_{-1,0}=(-1)^{[d]}4e_{-1,0}.
\end{eqnarray*}
According to $\LL_{-1}\otimes\QQ_0$, $\LL_0\otimes\QQ_{-1}$,
$\LL_1\otimes\QQ_{-2}$, $\QQ_0\otimes\LL_{-1}$,
$\QQ_{-2}\otimes\LL_1$, $\QQ_{-1}\otimes\LL_{0}$, one has
\begin{eqnarray*}
&&u_{-1,0}=u_{0,0}+(-1)^{[d]}2\mu_{-1,0},\
u_{-1,1}=(-1)^{[d]+1}4\mu_{-1,0},\
u_{-1,2}=(-1)^{[d]}\mu_{-1,0},\\
&&u^\dag_{-1,1}=u^\dag_{0,0}+(-1)^{[d]}2\mu_{-1,0},\
u^\dag_{-1,-1}=(-1)^{[d]}\mu_{-1,0},\
u^\dag_{-1,0}=(-1)^{[d]+1}4\mu_{-1,0}.
\end{eqnarray*}
Comparing the coefficients of $\GG_{-1}\otimes\HH_0$,
$\GG_0\otimes\HH_{-1}$, $\GG_1\otimes\HH_{-2}$,
$\HH_0\otimes\GG_{-1}$, $\HH_{-2}\otimes\GG_1$,
$\HH_{-1}\otimes\GG_{0}$, we obtain
$p^\dag_{-1,0}=p^\dag_{0,0}+(-1)^{[d]+1}\beta_{-1,0},\
p^\dag_{-1,1}=(-1)^{[d]}2\beta_{-1,0},\
p^\dag_{-1,2}=(-1)^{[d]+1}0.5\beta_{-1,0}$,
\begin{eqnarray*}
&&p_{-1,1}=p_{0,0}+(-1)^{[d]}\beta_{-1,0},\
p_{-1,-1}=(-1)^{[d]}0.5\beta_{-1,0},\
p_{-1,0}=(-1)^{[d]+1}2\beta_{-1,0}.
\end{eqnarray*}
Comparing the coefficients of $\GG_{-1}\otimes\QQ_0$,
$\GG_0\otimes\QQ_{-1}$, $\GG_1\otimes\QQ_{-2}$,
$\QQ_0\otimes\GG_{-1}$, $\QQ_{-2}\otimes\GG_1$,
$\QQ_{-1}\otimes\GG_{0}$, one has
$s_{-1,0}=s_{0,0}+(-1)^{[d]+1}e_{-1,0},\
s_{-1,1}=(-1)^{[d]}2e_{-1,0},\ s_{-1,2}=(-1)^{[d]+1}0.5e_{-1,0}$,
\begin{eqnarray*}
&&s^\dag_{-1,1}=s^\dag_{0,0}+(-1)^{[d]+1}e^\dag_{-1,0},\
s^\dag_{-1,-1}=(-1)^{[d]+1}0.5e^\dag_{-1,0},\
s^\dag_{-1,2}=(-1)^{[d]}2e^\dag_{-1,0}.
\end{eqnarray*}
Comparing the coefficients of $\LL_0\otimes\GG_{-1}$,
$\LL_{-2}\otimes\GG_1$, $\LL_{-1}\otimes\GG_{0}$,
$\GG_{-1}\otimes\LL_0$, $\GG_0\otimes\LL_{-1}$,
$\GG_{-2}\otimes\LL_1$, $\LL_{-1}\otimes\LL_0$,
$\LL_0\otimes\LL_{-1}$, $\LL_{-2}\otimes\LL_1$,
$\LL_1\otimes\LL_{-2}$, $\GG_{-1}\otimes\GG_0$,
$\GG_0\otimes\GG_{-1}$, $\GG_{-2}\otimes\GG_1$,
$\GG_1\otimes\GG_{-2}$, we have
\begin{eqnarray*}
x_{-1,0}+x_{1,0}=x_{-1,1}-x_{1,0}=x_{-1,-1}=x_{-1,2}=0\ \
\mbox{for}\,\,x=k,k^\dag,t,t^\dag.
\end{eqnarray*}
Comparing the coefficients of $\HH_{-1}\otimes\QQ_0$,
$\QQ_0\otimes\HH_{-1}$, $\HH_0\otimes\HH_{-1}$,
$\HH_{-1}\otimes\HH_0$, one has
\begin{eqnarray*}
q_{-1,0}=(-1)^{[d]}2\mu_{-1,0},\
q^\dag_{-1,1}=(-1)^{[d]}2\mu_{-1,0},\
w_{-1,1}=(-1)^{[d]+1}2e_{-1,0},\ w_{-1,0}=(-1)^{[d]}2e^\dag_{-1,0}.
\end{eqnarray*}

Applying $d_0$ to $[\LL_{2},\GG_{-1}]=3\GG_1$, we have $\LL_{2}\cdot
d_0(\GG_{-1})=(-1)^{[d]}\GG_{-1}\cdot d_0(\LL_{2})+3d_0(\GG_1)$.
Comparing the coefficients of $\LL_2\otimes\HH_{-1}$,
$\LL_0\otimes\HH_{1}$, $\HH_{-1}\otimes\LL_2$,
$\HH_{1}\otimes\LL_0$, $\LL_{2}\otimes\QQ_{-1}$,
$\LL_{0}\otimes\QQ_{1}$, $\GG_2\otimes\HH_{-1}$,
$\GG_0\otimes\HH_{1}$, $\LL_{1}\otimes\GG_{0}$,
$\GG_2\otimes\LL_{-1}$, $\LL_0\otimes\LL_{1}$, $\GG_0\otimes\GG_{1}$
and combine with (\ref{81031n01}), we have
$\theta_{\!\!-1,0}\!=\!\theta_{\!\!-2,0}\!=\!\theta_{2,0}\!=\!x_{\!1,0}\!=\!
y_{1,0}\!\!-\!\!y_{0,0}\!=\!z_{1,\!\!-1}\!\!-\!\!z_{0,0}\!=\!0$ for
$\theta\!=\!e^\dag,e,\mu,\beta,x\!=\!k,k^\dag,t,t^\dag,y\!=\!m,u,p^\dag,s,\
z\!=\!m^\dag,u^\dag,s^\dag,p$. Thus
$d_0(\LL_{-1})=d_0(\LL_{-2})=d_0(\LL_{2})=0$.

One can write
\begin{eqnarray*}
d_0(\HH_{n})\!\!\!&=\!\!\!&
\mbox{$\sum\limits_{i\in\Z}$}(\xi_{n,i}\LL_{i+n}\otimes\HH_{-i}+
\xi^\dag_{n,i}\HH_{i+n}\otimes\LL_{-i}+\eta_{n,i}\LL_{i+n}\otimes\GG_{-i}+
\eta^\dag_{n,i}\GG_{i+n}\otimes\LL_{-i}\\
&&+\varepsilon_{n,i}\LL_{i+n}\otimes\QQ_{-i}+\varepsilon^\dag_{n,i}\QQ_{i+n}\otimes\LL_{-i}
+\delta_{n,i}\HH_{i+n}\otimes\GG_{-i}+\delta^\dag_{n,i}\GG_{i+n}\otimes\HH_{-i}\\
&&+\psi_{n,i}\HH_{i+n}\otimes\QQ_{-i}+\psi^\dag_{n,i}\QQ_{i+n}\otimes\HH_{-i}
+\varphi_{n,i}\GG_{i+n}\otimes\QQ_{-i}+\varphi^\dag_{n,i}\QQ_{i+n}\otimes\GG_{-i}\\
&&+h_{n,i}\LL_{i+n}\otimes\LL_{-i}+l_{n,i}\HH_{i+n}\otimes\HH_{-i}+g_{n,i}\GG_{i+n}\otimes\GG_{-i}
+v_{n,i}\QQ_{i+n}\otimes\QQ_{-i}).
\end{eqnarray*}
Applying $d_0$ to $[\LL_{1},\HH_0]=0$, we have $\LL_{1}\cdot
d_0(\HH_0)=0$. Comparing the coefficients of
$\LL_{i+1}\otimes\HH_{-i}$, we have
$(i-1)\xi_{0,i}=(i+1)\xi_{0,i+1}\ \ \mbox{for}\ {i\in\Z}$, which
implies $\xi_{0,i}=0,\ \xi_{0,1}=-\xi_{0,0}$, $\forall\,\,i\neq0,1$.
Comparing the coefficients of $\LL_{i+1}\otimes\QQ_{-i}$,
$\GG_{i+1}\otimes\HH_{-i}$, $\GG_{i+1}\otimes\QQ_{-i}$,
$\HH_{i+1}\otimes\LL_{-i}$, $\QQ_{i+1}\otimes\LL_{-i}$,
$\HH_{i+1}\otimes\GG_{-i}$, $\QQ_{i+1}\otimes\GG_{-i}$, we have
\begin{eqnarray*}
x_{0,i}=y_{0,j}=x_{0,1}+x_{0,0}=y_{0,-1}+y_{0,0}=0,\
\,x=\varepsilon,\delta^\dag,\varphi,\
y=\xi^\dag,\varepsilon^\dag,\delta,\varphi^\dag,\forall\,\,
i\neq0,1,j\neq0,-1.
\end{eqnarray*}
Comparing the coefficients of $\LL_{i+1}\otimes\GG_{-i}$,
$\GG_{i+1}\otimes\LL_{-i}$, $\LL_{i+1}\otimes\LL_{-i}$,
$\GG_{i+1}\otimes\GG_{-i}$, we have
\begin{eqnarray*}
x_{0,i}=x_{0,-1}+\frac{1}{2}x_{0,0}= x_{0,1}-\frac{1}{2}x_{0,0}=0\
\mbox{for}\,x=\eta,\eta^\dag,h,g,\forall\,\,i\neq0,1,-1.
\end{eqnarray*}
Comparing the coefficients of $\HH_{i+1}\otimes\QQ_{-i}$,
$\QQ_{i+1}\otimes\HH_{-i}$, $\HH_{i+1}\otimes\HH_{-i}$,
$\QQ_{i+1}\otimes\QQ_{-i}$, one has
\begin{eqnarray*}
x_{0,i}=0\ \ \mbox{for}\ \,x=\psi,\psi^\dag,l,v, \forall\,\, i\neq0.
\end{eqnarray*}

Applying $d_0$ to $[\HH_{0},\GG_{0}]=\GG_0$, one has $\HH_{0}\cdot
d_0(\GG_{0})=(-1)^{[d]}\GG_{0}\cdot d_0(\HH_{0})+d_0(\GG_0)$.
Comparing the coefficients of $\LL_0\otimes\QQ_{0}$,
$\QQ_0\otimes\LL_{0}$, $\GG_{0}\otimes\HH_{0}$,
$\GG_{1}\otimes\HH_{-1}$, $\HH_{-1}\otimes\GG_1$,
$\GG_{0}\otimes\LL_{0}$, $\GG_{0}\otimes\GG_{0}$,
$\GG_{1}\otimes\GG_{-1}$, $\LL_{1}\otimes\LL_{-1}$,
$\LL_{-1}\otimes\LL_{1}$, $\GG_{0}\otimes\QQ_{0}$,
$\QQ_{0}\otimes\GG_{0}$, $\LL_{0}\otimes\HH_{0}$,
$\HH_{0}\otimes\LL_{0}$, one has
\begin{eqnarray*}
&&\!\!\!\!\!\!l_{0,0}\!=\!\varepsilon_{0,0}\!=\!\varepsilon^\dag_{0,0}=0,
u^\dag_{0,0}\!=\!-u_{0,0},\delta^\dag_{0,0}\!=\!\delta_{0,0},
\xi_{0,0}=\xi^\dag_{0,0}\!=\!-2\varphi_{0,0},
\varphi^\dag_{0,0}\!=\!-\varphi_{0,0},h_{0,0}\!=\!-4\varphi_{0,0},\\
&&\!\!\!\!\!\!v_{0,0}\!=\!(-1)^{[d]+1}u_{0,0},
\psi_{0,0}\!=\!(-1)^{[d]}s_{0,0},
\psi^\dag_{0,0}\!=\!(-1)^{[d]+1}s^\dag_{0,0},
m_{0,0}\!=\!2s^\dag_{0,0},\ m^\dag_{0,0}\!=\!-2s_{0,0}.
\end{eqnarray*}

Denoting $v_{1}=\QQ_{0}\otimes\HH_{0}$,
$v_{2}=\HH_{0}\otimes\QQ_{0}$ and observing the facts
\begin{eqnarray*}
\GG_{0}\cdot v_{1}=2\LL_{0}\otimes\HH_{0}+\QQ_{0}\otimes\GG_{0},\
\GG_{0}\cdot v_{2}=-\GG_{0}\otimes\QQ_{0}+2\HH_{0}\otimes\LL_{0},\
\LL_{\pm1}\cdot v_{i}=\LL_{\pm2}\cdot v_{i}=0,
\end{eqnarray*}
for $i=1,2$, replacing $d_0$ by $d_0-\kappa_{\rm inn}$ (where
$\kappa$ is some linear combination of $v_{1}$, $v_{2}$), one has
\begin{eqnarray*}
d_0(\GG_{0})\!\!\!&=&\!\!\!p^\dag_{0,0}\GG_{0}\otimes\HH_{0}+p_{0,0}\HH_{0}\otimes\GG_{0}
+u_{0,0}\LL_{0}\otimes\QQ_{0}-u_{0,0}\QQ_{0}\otimes\LL_{0},\\
d_0(\GG_{1})\!\!\!&=&\!\!\!-u_{0,0}\QQ_{0}\otimes\LL_{1}+u_{0,0}\LL_{1}\otimes\QQ_{0}
+p_{0,0}\HH_{0}\otimes\GG_{1}+p^\dag_{0,0}\GG_{1}\otimes\HH_{0},\\
d_0(\GG_{-1})\!\!\!&=&\!\!\!u_{0,0}\LL_{-1}\otimes\QQ_{0}
-u_{0,0}\QQ_{0}\otimes\LL_{-1}+p^\dag_{0,0}\GG_{-1}\otimes\HH_{0}
+p_{0,0}\HH_{0}\otimes\GG_{-1}.
\end{eqnarray*}

Comparing the coefficients of $\GG_0\otimes\HH_{-1}$,
$\GG_{-2}\otimes\GG_{1}$ in $\HH_{0}\cdot
d_0(\GG_{-1})=(-1)^{[d]}\GG_{-1}\cdot d_0(\HH_{0})+d_0(\GG_{-1})$,
we have $\varphi_{0,0}=\delta_{0,0}=0$. Comparing the coefficients
of $\LL_{-1}\otimes\GG_{0}$, $\GG_{-1}\otimes\LL_{0}$,
$\GG_{-1}\otimes\GG_{0}$ in $\LL_{-1}\cdot d_0(\HH_{0})=\HH_{0}\cdot
d_0(\LL_{-1})$, we have $\eta_{0,0}=\eta^\dag_{0,0}=g_{0,0}=0$. Thus
one can suppose
\begin{eqnarray*}
d_0(\HH_{0})=(-1)^{[d]+1}u_{0,0}\QQ_{0}\otimes\QQ_{0}.
\end{eqnarray*}
Write
\begin{eqnarray*}
d_0(\QQ_{n})\!\!\!&=\!\!\!&
\mbox{$\sum\limits_{i\in\Z}$}(A_{n,i}\LL_{i+n}\otimes\HH_{-i}+
A^\dag_{n,i}\HH_{i+n}\otimes\LL_{-i}+B_{n,i}\LL_{i+n}\otimes\GG_{-i}+
B^\dag_{n,i}\GG_{i+n}\otimes\LL_{-i}\\
&&+D_{n,i}\LL_{i+n}\otimes\QQ_{-i}+D^\dag_{n,i}\QQ_{i+n}\otimes\LL_{-i}
+E_{n,i}\HH_{i+n}\otimes\GG_{-i}+E^\dag_{n,i}\GG_{i+n}\otimes\HH_{-i}\\
&&+F_{n,i}\HH_{i+n}\otimes\QQ_{-i}+F^\dag_{n,i}\QQ_{i+n}\otimes\HH_{-i}
+R_{n,i}\GG_{i+n}\otimes\QQ_{-i}+R^\dag_{n,i}\QQ_{i+n}\otimes\GG_{-i}\\
&&+X_{n,i}\LL_{i+n}\otimes\LL_{-i}+Y_{n,i}\HH_{i+n}\otimes\HH_{-i}+M_{n,i}\GG_{i+n}\otimes\GG_{-i}
+N_{n,i}\QQ_{i+n}\otimes\QQ_{-i}).
\end{eqnarray*}

Comparing the coefficients of $\GG_0\otimes\HH_{0}$,
$\HH_0\otimes\GG_{0}$, $\LL_{0}\otimes\GG_{0}$,
$\GG_{0}\otimes\LL_{0}$, $\GG_0\otimes\GG_{0}$,
$\QQ_0\otimes\QQ_{0}$, $\LL_0\otimes\HH_{0}$, $\HH_0\otimes\LL_{0}$,
$\QQ_0\otimes\GG_{0}$, $\GG_0\otimes\QQ_{0}$, $\LL_0\otimes\LL_{0}$,
$\HH_0\otimes\HH_{0}$ in $\HH_{0}\cdot
d_0(\QQ_{0})+d_0(\QQ_0)=(-1)^{[d]}\QQ_{0}\cdot d_0(\HH_{0})$, we
obtain
$E^\dag_{0,0}\!=\!E_{0,0}=B_{0,0}\!=\!B^\dag_{0,0}\!=\!M_{0,0}\!=\!X_{0,0}\!=\!Y_{0,0}\!=\!A_{0,0}\!=\!A^\dag_{0,0}
\!=\!R^\dag_{0,0}\!=\!R_{0,0}\!=\!N_{0,0}\!=\!0$. Thus we can
rewrite $d_0(\QQ_{0})$ as
\begin{eqnarray*}
D_{0,0}\LL_{0}\otimes\QQ_{0}-D_{0,0}\LL_{1}\otimes\QQ_{-1}
+D^\dag_{0,0}\QQ_{0}\otimes\LL_{0}-D^\dag_{0,0}\QQ_{-1}\otimes\LL_{1}+F_{0,0}\HH_{0}\otimes\QQ_{0}+F^\dag_{0,0}\QQ_{0}\otimes\HH_{0}.
\end{eqnarray*}
Applying $d_0$ to $[\GG_{0},\QQ_{0}]=2\LL_{0}$, we have
$\GG_{0}\cdot d_0(\QQ_{0})+\QQ_{0}\cdot d_0(\GG_{0})=0$. Comparing
the coefficients of $\GG_{0}\otimes\QQ_{0}$,
$\QQ_{0}\otimes\GG_{0}$, $\HH_{0}\otimes\LL_{0}$,
$\LL_{0}\otimes\HH_{0}$, $\GG_{1}\otimes\QQ_{-1}$,
$\HH_{-1}\otimes\LL_{1}$, we have $F_{0,0}=F^\dag_{0,0}=-p_{0,0},\
p^\dag_{0,0}=p_{0,0}, D_{0,0}=D^\dag_{0,0}=0$. Thus one can suppose
\begin{eqnarray*}
d_0(\QQ_{0})\!\!\!&=&\!\!\!
-p_{0,0}\HH_{0}\otimes\QQ_{0}-p_{0,0}\QQ_{0}\otimes\HH_{0},\\
d_0(\GG_{0})\!\!\!&=&\!\!\!
p_{0,0}\GG_{0}\otimes\HH_{0}+p_{0,0}\HH_{0}\otimes\GG_{0}
+u_{0,0}\LL_{0}\otimes\QQ_{0}-u_{0,0}\QQ_{0}\otimes\LL_{0}.
\end{eqnarray*}
Denoting $x=\HH_{0}\otimes\HH_{0}$, observing $\GG_{0}\cdot x
=-\GG_{0}\otimes\HH_{0}-\HH_{0}\otimes\GG_{0}$ and replacing $d_0$
by $d_0+p_{0,0}x_{\rm inn}$ (noting that this replacement does not
affect $d_0(\LL_{\pm1}),\ d_0(\LL_{\pm2})$), one can suppose
\begin{eqnarray*}
&&\!\!\!\!\!\!\!\!\!\!\!\!d_0(\QQ_{0})=0,\
d_0(\HH_{0})=(-1)^{[d]+1}u_{0,0}\QQ_{0}\otimes\QQ_{0},\
d_0(\GG_{0})=u_{0,0}\LL_{0}\otimes\QQ_{0}
-u_{0,0}\QQ_{0}\otimes\LL_{0},\\
&&\!\!\!\!\!\!\!\!\!\!\!\!d_0(\GG_{1})=u_{0,0}\LL_{1}\otimes\QQ_{0}
-u_{0,0}\QQ_{0}\otimes\LL_{1},\ d_0(\GG_{-1})=
u_{0,0}\LL_{-1}\otimes\QQ_{0}-u_{0,0}\QQ_{0}\otimes\LL_{-1}.
\end{eqnarray*}
Applying $d_0$ to $[\LL_{-1},\QQ_{1}]=-\QQ_{0}$, we have
$\LL_{-1}\cdot d_0(\QQ_{1})+d_0(\QQ_{0})\!=\!(-1)^{[d]}\QQ_{1}\cdot
d_0(\LL_{-1})$. Comparing the coefficients of
$\LL_{i}\otimes\HH_{-i}$, we obtain $(i+2)A_{1,i}=(i-1)A_{1,i-1}$,
$\forall\,\, i\neq0,1,2$ and then
$A_{1,i}=A_{1,-1}+2A_{1,0}=A_{1,-2}-A_{1,0}=0$, $\forall\,\,
i\neq-1,-2,0$. Similarly, comparing the coefficients of
$\LL_{i}\otimes\QQ_{-i}$, $\GG_{i}\otimes\HH_{-i}$,
$\GG_{i}\otimes\QQ_{-i}$, we obtain
\begin{eqnarray*}
x_{1,i}=x_{1,-1}+2x_{1,0}=x_{1,-2}-x_{1,0}=0\
\mbox{for}\,x=D,E^\dag,R, \forall\,\, i\neq-1,-2,0.
\end{eqnarray*}
Comparing the coefficients of $\HH_{i}\otimes\LL_{-i}$, one has
$(i+1)A^\dag_{1,i}=(i-2)A^\dag_{1,i-1}$, $\forall\,\, i\neq-1,-2,0$,
which implies.
$A^\dag_{1,i}=2A^\dag_{1,-1}+A^\dag_{1,0}=2A^\dag_{1,1}+A^\dag_{1,0}=0$,
$\forall\,\, i\neq-1,1,0$. Similarly comparing the coefficients of
$\QQ_{i}\otimes\LL_{-i}$, $\HH_{i}\otimes\GG_{-i}$,
$\QQ_{i}\otimes\GG_{-i}$, we have
\begin{eqnarray*}
x_{1,i}=2x_{1,-1}+x_{1,0}=2x_{1,1}+x_{1,0}=0,\
\mbox{for}\,x=D^\dag,E,R^\dag, \forall\,\, i\neq-1,1,0.
\end{eqnarray*}
Comparing the coefficients of $\LL_{i}\otimes\GG_{-i}$, one has
$(i+2)B_{1,i}=(i-2)B_{1,i-1}$ and
\begin{eqnarray*}
B_{1,i}=B_{1,-1}+B_{1,0}=3B_{1,1}+B_{1,0}=3B_{1,-2}-B_{1,0}=0,\
\forall\,\, i\neq-1,-2,0,1.
\end{eqnarray*}
Comparing the coefficients of $\GG_{i}\otimes\LL_{-i}$,
$\LL_{i}\otimes\LL_{-i}$, $\GG_{i}\otimes\GG_{-i}$, one has
\begin{eqnarray*}
x_{1,i}=x_{1,-1}+x_{1,0}=3x_{1,1}+x_{1,0}=3x_{1,-2}-x_{1,0}=0,\
\mbox{for}\,x=X,M,B^\dag, \forall\,\, i\neq-1,-2,0,1.
\end{eqnarray*}
Comparing the coefficients of $\HH_{i}\otimes\QQ_{-i}$, we obtain
$(i+1)F_{1,i}=(i-1)F_{1,i-1}$, $\forall\,\, i\neq-1,-2,0,1,2$ and
$F_{1,i}=F_{1,-1}+F_{1,0}=0$, $\forall\,\, i\neq-1,0$. Similarly,
comparing the coefficients of $\QQ_{i}\otimes\HH_{-i}$,
$\HH_{i}\otimes\HH_{-i}$, $\QQ_{i}\otimes\QQ_{-i}$, we have
$x_{1,i}=x_{1,-1}+x_{1,0}=0$ for $x=F^\dag,Y,N$, $\forall\,\,
i\neq-1,0$. Applying $d_0$ to $[\HH_{0},\QQ_{1}]=-\QQ_{1}$, we have
$\HH_{0}\cdot d_0(\QQ_{1})+d_0(\QQ_{1})=(-1)^{[d]}\QQ_{1}\cdot
d_0(\HH_{0})$. Comparing the coefficients of
$\GG_{1}\otimes\HH_{0}$, $\HH_{0}\otimes\GG_{1}$,\
$\LL_{1}\otimes\GG_{0}$,\ $\GG_{1}\otimes\LL_{0}$,\
$\GG_{1}\otimes\GG_{0}$,\ $\GG_{1}\otimes\QQ_{0}$,\
$\QQ_{0}\otimes\QQ_{1}$,\ $\LL_{1}\otimes\HH_{0}$,\
$\GG_{1}\otimes\QQ_{0}$,\ $\HH_{1}\otimes\LL_{0}$,\
$\QQ_{1}\otimes\GG_{0}$, one has
\begin{eqnarray*}
E^\dag_{1,0}=E_{1,0}=B_{1,0}=B^\dag_{1,0}=M_{1,0}=A^\dag_{1,0}=X_{1,0}=A_{1,0}=R_{1,0}=R^\dag_{1,0}=Y_{1,0}=N_{1,0}=0.
\end{eqnarray*}
Applying $d_0$ to $[\GG_{-1},\QQ_{1}]=2\LL_{0}-2\HH_{0}$, we have
$(-1)^{[d]}\GG_{-1}\cdot d_0(\QQ_{1})+(-1)^{[d]}\QQ_{1}\cdot
d_0(\GG_{-1})+2d_0(\HH_{0})=0$. Comparing the coefficients of
$\GG_{0}\otimes\QQ_{0}$, $\QQ_{1}\otimes\GG_{-1}$, one has
$2D_{1,0}+F_{1,0}=D^\dag_{1,0}+F^\dag_{1,0}=0$. Applying $d_0$ to
$[\LL_{1},\QQ_{-1}]=\QQ_{0}$, we have
\begin{eqnarray}\label{081124a01}
\LL_{1}\cdot d_0(\QQ_{-1})\!\!\!&=\!\!\!& d_0(\QQ_{0}).
\end{eqnarray}
Comparing the coefficients of $\LL_{i}\otimes\HH_{-i}$ in
(\ref{081124a01}), we obtain $(i-2)A_{-1,i}=(i+1)A_{-1,i+1}$,
$\forall\,\, i\neq0,1$ and
$A_{-1,i}=A_{-1,1}+2A_{-1,0}=A_{-1,2}-A_{-1,0}=0$, $\forall\,\,i
\neq1,0,2$. Comparing the coefficients of $\LL_{i}\otimes\QQ_{-i}$,
$\GG_{i}\otimes\HH_{-i}$, $\GG_{i}\otimes\QQ_{-i}$, we have
\begin{eqnarray*}
&&x_{-1,i}=x_{-1,1}+2x_{-1,0}=x_{-1,2}-x_{-1,0}=0\
\mbox{for}\,x=D,E^\dag,R,\ \forall\,\,i \neq1,0,2.
\end{eqnarray*}
Comparing the coefficients of $\HH_{i}\otimes\LL_{-i}$ in
(\ref{081124a01}), one has
$A^\dag_{-1,1}=A^\dag_{-1,-1}=-\frac{1}{2}A^\dag_{-1,0}$,
$\forall\,\, i\neq1,-1,0$. Comparing the coefficients of
$\QQ_{i}\otimes\LL_{-i}$, $\HH_{i}\otimes\GG_{-i}$,
$\QQ_{i}\otimes\GG_{-i}$, we have
\begin{eqnarray*}
&&x_{-1,i}=0,\ x_{-1,1}=x_{-1,-1}=-\frac{1}{2}x_{-1,0}\
\mbox{for}\,x=D^\dag,E,R^\dag, \forall\,\, i\neq1,-1,0.
\end{eqnarray*}
Comparing the coefficients of $\LL_{i}\otimes\GG_{-i}$, we obtain
$(i-2)B_{-1,i}=(i+2)B_{-1,i+1}$ for $i\in\Z$ and
\begin{eqnarray*}
B_{-1,i}=B_{-1,1}+B_{-1,0}=
3B_{-1,-1}+B_{-1,0}=3B_{-1,-2}-B_{-1,0}=0,\ \forall\,\,
i\neq1,-1,-2,0.
\end{eqnarray*}
Comparing the coefficients of $\GG_{i}\otimes\LL_{-i}$,
$\LL_{i}\otimes\LL_{-i}$, $\GG_{i}\otimes\GG_{-i}$, one has
$x_{-1,i}=x_{-1,1}+x_{-1,0}=
3x_{-1,-1}+x_{-1,0}=3x_{-1,-2}-x_{-1,0}=0$ for $x=X,M,B^\dag$,
$\forall\,\, i\neq1,-1,-2,0$. Comparing the coefficients of
$\HH_{i}\otimes\QQ_{-i}$ in (\ref{081124a01}), we have
$(i-1)F_{-1,i}=(i+1)F_{-1,i+1}$, $\forall\,\, i\neq0$ and
$F_{-1,i}=F_{-1,1}+F_{-1,0}=0$, $\forall\,\, i\neq1,0$. Comparing
the coefficients of $\QQ_{i}\otimes\HH_{-i}$,
$\HH_{i}\otimes\HH_{-i}$, $\QQ_{i}\otimes\QQ_{-i}$, one has
$x_{-1,i}=x_{-1,1}+x_{-1,0}=0$ for $x=F^\dag,Y,N$, $\forall\,\,
i\neq1,0$.

Applying $d_0$ to $[\HH_{0},\QQ_{-1}]=\QQ_{-1}$, we have
$\HH_{0}\cdot d_0(\QQ_{-1})=(-1)^{[d]}\QQ_{-1}\cdot
d_0(\HH_{0})-d_0(\QQ_{-1})$. Comparing the coefficients of
$\GG_{-1}\otimes\HH_{0}$, $\LL_{-1}\otimes\GG_{0}$,
$\GG_{-1}\otimes\LL_{0}$, $\GG_{-1}\otimes\GG_{0}$,
$\QQ_{0}\otimes\QQ_{-1}$, $\LL_{-1}\otimes\HH_{0}$,
$\HH_{-1}\otimes\HH_{0}$, $\GG_{-1}\otimes\QQ_{0}$,
$\QQ_{-1}\otimes\GG_{0}$, $\HH_{-1}\otimes\LL_{0}$,
$\LL_{-1}\otimes\LL_{0}$, one has
$E_{-1,0}=E^\dag_{-1,0}=B_{-1,0}=B^\dag_{-1,0}=M_{-1,0}
=R_{-1,0}=R^\dag_{-1,0}=A^\dag_{-1,0}=X_{-1,0}=N_{-1,0}=A_{-1,0}=Y_{-1,0}=0$.
Applying $d_0$ to $[\LL_{2},\QQ_{-1}]=\QQ_{1}$, we have
$\LL_{2}\cdot d_0(\QQ_{-1})=(-1)^{[d]}\QQ_{-1}\cdot
d_0(\LL_{2})+d_0(\QQ_{1})$. Comparing the coefficients of
$\LL_{2}\otimes\QQ_{-1}$, $\LL_{1}\otimes\QQ_{2}$,
$\QQ_{-2}\otimes\LL_{3}$, $\QQ_{2}\otimes\LL_{-1}$,
$\HH_{1}\otimes\QQ_{0}$, $\QQ_{1}\otimes\HH_{0}$, one has
\begin{eqnarray*}
D_{-1,0}=D_{1,0}=D^\dag_{-1,0}=D^\dag_{1,0}=F_{-1,0}-F_{1,0}=
F^\dag_{-1,0}-F^\dag_{1,0}=0.
\end{eqnarray*}
Then $F^\dag_{1,0}= F_{1,0}=F^\dag_{-1,0}=F_{-1,0}=0$ and
$d_0(\QQ_{1})=d_0(\QQ_{-1})=0$.

Applying $d_0$ to $[\LL_{1},\HH_{-1}]=\HH_{0}$, we have
$\LL_{1}\cdot d_0(\HH_{-1})=d_0(\HH_{0})$. Comparing the
coefficients of $\LL_{i}\otimes\HH_{-i}$, we have
$(i-2)\xi_{-1,i}=(i+1)\xi_{-1,i+1}\ \ \mbox{for}\ {i\in\Z}$. we
obtain $\xi_{-1,i}=\xi_{-1,1}+2\xi_{-1,0}=\xi_{-1,2}-\xi_{-1,0}=0,\
\forall\,\,i\neq0,1,2$. Similarly comparing the coefficients of
$\LL_{i}\otimes\QQ_{-i}$, $\GG_{i}\otimes\QQ_{-i}$,
$\LL_{i}\otimes\QQ_{-i}$, we have
$x_{-1,i}=x_{-1,1}+2x_{-1,0}=x_{-1,2}-x_{-1,0}=0$ for
$x=\delta^\dag,\varphi,\varepsilon, \forall\,\,i\neq0,1,2$.
Comparing the coefficients of $\HH_{i}\otimes\LL_{-i}$, we have
$(i-1)\xi^\dag_{-1,i}=(i+2)\xi^\dag_{-1,i+1}\ \ \mbox{for}\
{i\in\Z}$ and
\begin{eqnarray*}
\xi^\dag_{-1,i}=0,\
\xi^\dag_{-1,1}=\xi^\dag_{-1,-1}=-\frac{1}{2}\xi^\dag_{-1,0},\
\forall\,\,i\neq1,-1,0.
\end{eqnarray*}
Comparing the coefficients of $\HH_{i}\otimes\GG_{-i}$,
$\QQ_{i}\otimes\GG_{-i}$, $\QQ_{i}\otimes\LL_{-i}$, we have
$x_{-1,i}=0,\ x_{-1,1}=x_{-1,-1}=-\frac{1}{2}x_{-1,0}$ for
$x=\delta,\varphi^\dag,\varepsilon^\dag$, $\forall\,\, i\neq1,-1,0$.
Comparing the coefficients of $\LL_{i}\otimes\GG_{-i}$, we have
$(i-2)\eta_{-1,i}=(i+2)\eta_{-1,i+1}\ \ \mbox{for}\ {i\in\Z}$ and
\begin{eqnarray*}
\eta_{-1,i}=\eta_{-1,1}+\eta_{-1,0}=3\eta_{-1,-1}+\eta_{-1,0}=
3\eta_{-1,2}-\eta_{-1,0}=0,\  \forall\,\,i\neq1,-1,2,0.
\end{eqnarray*}
Comparing the coefficients of $\GG_{i}\otimes\LL_{-i}$,\
$\LL_{i}\otimes\LL_{-i}$, $\GG_{i}\otimes\GG_{-i}$, we have
\begin{eqnarray*}
x_{-1,i}=x_{-1,1}+x_{-1,0}=3x_{-1,-1}+x_{-1,0}=
3x_{-1,2}-x_{-1,0}=0,\ x=\eta^\dag,h,g,\  \forall\,\, i\neq0,\pm1,2.
\end{eqnarray*}
Comparing the coefficients of $\HH_{i}\otimes\QQ_{-i}$, we obtain
$(i-1)\psi_{-1,i}=(i+1)\psi_{-1,i+1}\ \mbox{for}\ i\in\Z^*$ and
$\psi_{-1,i}=\psi_{-1,1}+\psi_{-1,0}=0,\ \forall\,\,i\neq1,0$.
Comparing the coefficients of $\QQ_{i}\otimes\HH_{-i}$,
$\HH_{i}\otimes\HH_{-i}$, we have $x_{-1,i}=x_{-1,1}+x_{-1,0}=0$ for
$x=\psi^\dag,l$, $\forall\,\, i\neq1,0$. According to
$\QQ_{i}\otimes\QQ_{-i}$, we obtain $(i-1)v_{-1,i}=(i+1)v_{-1,i+1}\
\ \mbox{for}\ i\neq0$ and
$v_{-1,i}=v_{-1,1}+v_{-1,0}-(-1)^{[d]+1}u_{0,0}=0$,
$\forall\,\,i\neq1,0$.

Applying $d_0$ to $[\HH_{-1},\GG_{1}]=\GG_{0}$, we have
$\HH_{-1}\cdot d_0(\GG_{1})=(-1)^{[d]}\GG_{1}\cdot
d_0(\HH_{-1})+d_0(\GG_{0})$. Comparing the coefficients of
$\GG_{1}\otimes\QQ_{-1}$, $\LL_{0}\otimes\LL_{0}$, one has
$2\varepsilon_{-1,0}-\psi_{-1,0}=2\varepsilon_{-1,0}-\varepsilon^\dag_{-1,0}=0$.
Applying $d_0$ to $[\HH_{-1},\QQ_{1}]=-\QQ_{0}$, we have
$\HH_{-1}\cdot d_0(\QQ_{1})+d_0(\QQ_{0})=(-1)^{[d]}\QQ_{1}\cdot
d_0(\HH_{-1})$. Comparing the coefficients of
$\HH_{1}\otimes\QQ_{-1}$, $\QQ_{1}\otimes\HH_{-1}$,
$\QQ_{0}\otimes\QQ_{0}$, $\QQ_{2}\otimes\QQ_{-2}$,
$\QQ_{2}\otimes\HH_{-2}$, $\HH_{-2}\otimes\QQ_{2}$,
$\HH_{1}\otimes\HH_{-1}$, $\HH_{-1}\otimes\HH_{1}$,
$\QQ_{-2}\otimes\LL_{2}$, $\LL_{-2}\otimes\LL_{2}$,
$\LL_{2}\otimes\LL_{-2}$, $\QQ_{-1}\otimes\QQ_{1}$, we have
$\varepsilon_{-1,0}\!=\!l_{-1,0}\!=\!\xi_{-1,0}\!
=\!\xi^\dag_{-1,0}\!=\!\delta_{-1,0}\!
=\!\delta^\dag_{-1,0}\!=\!\eta_{-1,0}\!=\!\eta^\dag_{-1,0}\!
=\!\varphi^\dag_{-1,0}\!=\!h_{-1,0}\!=\!g_{-1,0}\!
=\!\varphi_{-1,0}\!=\!\varepsilon^\dag_{-1,0}\!=\!\psi_{-1,0}\!=\!\psi^\dag_{-1,0}\!=\!0$.
Applying $d_0$ to $[\HH_{-1},\HH_{0}]=0$, we have $\HH_{-1}\cdot
d_0(\HH_{0})=\HH_{0}\cdot d_0(\HH_{-1})$. Comparing the coefficients
of $\QQ_{0}\otimes\QQ_{-1}$, one has $u_{0,0}=0$. Then
$d_0(\GG_{\pm1})\!=\!d_0(\GG_{0})\!=\!d_0(\HH_{0})\!=\!d_0(\HH_{-1})\!=\!0$.
\begin{clai}\adddot\rm\label{add-clai+}
The sum in (\ref{summable}) is finite.
\end{clai}
By now, we have completed the proof of Proposition \ref{prop3.2}.

The following lemma can be obtained by using the similar techniques
developed in \cite{YS}.
\begin{lemm}\adddot\label{lemma4}
If $x\cdot r\in\Im(1\otimes1-\tau)$ for some $r\in V$
and any $x\in \TT$, then $r\in\Im(1\otimes1-\tau)$.
\end{lemm}

\ni{\it Proof of Theorem \ref{main}.} Suppose $(\TT , [\cdot,\cdot],
\triangle)$ is a Lie super-bialgebra structure on $\TT$. Then
$\Delta=\Delta_{r}$ for some $r\in (\TT \otimes \TT) ^{\bar{0}}$ by
Proposition \ref{prop3.2}. Then $\Im\,\Delta \subset \Im(1\otimes1-
\tau)$. Thus by Lemma \ref{lemma4}, we have
$r\in\Im(1\otimes1-\tau)$. Then Lemma \ref{lemm3.1}(iii) shows that
$c (r)=0$. Hence $(\TT , [\cdot,\cdot], \triangle)$ is a triangular
coboundary Lie super-bialgebra.

\end{document}